\documentclass[11pt]{article}
\usepackage{amsmath, amsthm}
\usepackage{amsfonts}
\usepackage{amssymb}
\usepackage{euscript}
\usepackage[square, comma, sort&compress, numbers]{natbib}

\let\oldref=\ref
\def\ref#1{(\oldref{#1})}

\newcommand\blfootnote[1]{%
  \begingroup
  \renewcommand\thefootnote{}\footnote{#1}%
  \addtocounter{footnote}{-1}%
  \endgroup
}

\theoremstyle{plain}
\newtheorem{thm}{Theorem}[section]

\newtheorem{cor}[thm]{Corollary}
\newtheorem{lem}[thm]{Lemma}
\newtheorem{prop}[thm]{Proposition}
\theoremstyle{definition}
\newtheorem{ex}[thm]{Example}
\newtheorem{defn}[thm]{Definition}
\newtheorem{nota}[thm]{Notation}
\newtheorem{rem}[thm]{Remark}

\newcommand{\Z}{\mathrm{Z}}
\newcommand{\J}{\mathrm{J}}
\newcommand{\U}{\mathrm{U}}
\newcommand{\Ann}{\mathrm{Ann}}

\newcommand{\fm}{\mathfrak{M}}

\newcommand{\cp}{\EuScript{P}}

\newcommand{\n}{\mathbb{N}}

\newcommand{\pre}{pr\'e\-sim\-pli\-fi\-able}

\newcommand{\ifof}{if and only if }
\newcommand{\irr}{ $\{irreducible, strongly irreducible, very strongly irreducible\}$ }
\newcommand{\prim}{ $\{primitive, strongly primitive, very strongly primitive\}$ }

\newcommand{\se}{\subseteq}

\newcommand{\sm}{\setminus}

\newcommand{\give}{$\Rightarrow$}
\newcommand{\rgive}{$\Leftarrow$}

\newcommand{\f}[2]{\frac{#1}{#2}}

\newcommand{\tohi}{\emptyset}

\notag
                                                  
\begin{document}

\title{Factorizations in Modules and Splitting Multiplicatively Closed Subsets}
\author{Ashkan Nikseresht \\
\it\small Department of Mathematics, Institute for Advanced Studies in Basic Sciences,\\
\it\small P.O. Box 45195-1159, Zanjan, Iran\\
\small E-mail:  ashkan\_nikseresht@yahoo.com%
}
\date{}

\maketitle

\begin{abstract}
We introduce the concept of multiplicatively closed subsets of a commutative ring $R$ which split an $R$-module
$M$ and study factorization properties of elements of $M$ with respect to such a set. Also we demonstrate  how
one can utilize this concept to investigate factorization properties of $R$ and deduce some Nagata type theorems
relating factorization properties of $R$ to those of its localizations, when $R$ is an integral domain.
\end{abstract}


                                         \section{Introduction}
Throughout this paper all rings are commutative with identity and all modules are unitary. We assume that all
modules are nonzero. Also $R$ denotes a ring and $M$ is an $R$-module.

Theory of factorization in commutative rings which has a long history (see for example \cite{Nagata}), still gets
a lot of attention from various researchers. To see some recent papers on this subject, the reader is referred to
\cite{bfr, andnew1, QFD=UFD, baeth2, chang, chun, fooua, Zn, we, smcs, thez, FFR}. In \cite{r1, r2}, D. D.
Anderson and S. Valdes-Leon generalized the theory of factorization in integral domains to commutative rings with
zero divisors and to modules as well. These concepts are further studied in \cite{bfr, chun, FFR, r3, we}.

One of the longstanding questions in this subject is ``what is the relationship between factorization properties
of $R$ and those of its localizations?'', especially when $R$ is a domain (see for example \cite{d2, r3, Nagata}).
In particular, many have tried to give conditions under which, if $R_S$ is a UFD (or has other types of
factorization properties), then $R$ is so, where $S$ is a multiplicatively closed subset of $R$. For
example \cite[Corollary 8.32]{lar}, says that if $R$ is a Krull domain and $S$ is generated by a set of primes and
$R_S$ is a UFD, then $R$ is a UFD. This type of results, are called Nagata type theorems due to a theorem of
Nagata in \cite{Nagata}. One can find some other similar results and a brief review of this subject in
\cite[Section 3]{d2}.

On the other hand, in \cite{smcs} the concept of factorization with respect to a saturated multiplicatively closed
subset (also called a divisor closed multiplicative submonoid) of $R$  is introduced. If we apply Theorem 3.9 of
that paper with $S'=M=R$ and assume that $R$ is an integral domain, then we get a Nagata type result which states
that if $R_S$ is a bounded factorization domain and $R$ as an $R$-module is an $S$-bounded factorization module,
then $R$ is a bounded factorization domain (for exact definitions, see Section 2). It is still unknown whether the
similar result holds for other factorization properties such as unique or finite factorization (see \cite[Question
3.11]{smcs}).

The main aim of  this research is to find partial answers to the above question and utilize them to find relations
between factorization properties of $R$ and $R_S$, especially when $R$ is a domain. For this, we generalize the
concept of an splitting multiplicatively closed subset in \cite{d2} which is of key importance in the results of
that paper. Interestingly, we find out that this concept is equivalent to another one which is completely stated
in terms of factorization properties with respect to a saturated multiplicatively closed subset.

In Section 2, we briefly review the concepts of factorization theory with respect to a saturated multiplicatively
closed subset. Then in Section 3, we state the definition of a multiplicatively closed subset which splits $M$ and
study basic properties of such sets. In Section 4, we present our main results, which state how factorization
properties of $M$ and $M_S$ are related when $S$ splits $M$. Finally, in Section 5, we present an example in which
$M=R$ is an integral domain to show how our result could be applied in order to study factorization properties of
integral domains.

In the following, by $\U(R)$ and $\J(R)$ we mean the set of units and Jacobson radical of $R$, respectively.
Furthermore, $\Z(N)$, where $N\se M$, means the set of zero divisors of $N$, that is, $\{r\in R|\exists 0\neq m\in
N: \quad rm=\ 0\}$. In addition, $\Ann(N)$ (resp. $\Ann_M(r)$) denotes the annihilator in $R$ of $N\se M$ (resp.
in $M$ of $r\in R$). Any other undefined notation is as in \cite{aty}.

               \section{A brief review of factorization with respect to a saturated multiplicatively closed subset}

In this section we recall the main concepts of factorization with respect to a saturated multiplicatively closed
subset of $R$ which is needed in this paper. For more details and several examples the reader is referred to
\cite{smcs}. In the sequel, $S$ always
denotes a saturated multiplicatively closed subset of $R$ (we let $S$ to contain 0, which means $S=R$). 

We say that two elements, $m$ and $n$ of $M$, are \emph{$S$-associates} and write $m\sim^S n$, if there exist
$s,s'\in S$ such that $m=sn$ and $n=s'm$. They are called \emph{$S$-strong associates}, if $m=un$ for some $u\in
\U(R)\cap S$ and we denote it by $m\approx^S n$. Also we call them \emph{$S$-very strong associates}, denoted by
$m\cong^S n$, when $m\sim^S n$ and either $m=n=0$ or from $m=sn$ for some $s\in S$ it follows that $s\in \U(R)$.

In the case that $S=R$, we drop the $S$ prefixes. In this case, our notations coincide with that of \cite{r1, r2}.
An $m\in M$ is called \emph{$S$-primitive} (resp. \emph{$S$-strongly primitive}, \emph{$S$-very strongly
primitive}), when $m=sn$ for some $s\in S,\;n\in M$ implies $n\sim^S m$ (resp. $n\approx^S m$, $n\cong^S m$). A
nonunit element $a\in R$ is called \emph{irreducible} (resp. \emph{strongly irreducible}, \emph{very strongly
irreducible}) if $a=bc$ for some $b,c\in R$, implies $a\sim b$ or $a\sim c$ (resp. $a\approx b$ or $a\approx c$,
$a\cong b$ or $a\cong c$). Note that here by being associates in $R$, we mean being associates in $R$ as an
$R$-module.


By an \emph{$S$-factorization} of $m\in M$ with length $k$, we mean an equation $m=s_1\cdots s_kn$ where $s_i$'s
are nonunits in $S$, $k\in\n\cup\{0\}$ and $n\in M$. If moreover, for some $\alpha\in \irr$ and  $\beta\in \prim$
$s_i$'s are $\alpha$ and $n$ is $S$-$\beta$, we call this an \emph{($\alpha,\,\beta$)-$S$-factorization}. If every
nonzero element of $M$ has an ($\alpha,\, \beta$)-$S$-factorization, we say that $M$ is
\emph{($\alpha,\,\beta$)-$S$-atomic}.

By an \emph{$S$-atomic factorization} we mean an (irreducible, primitive)-$S$-factorization and by an
\emph{$S$-atomic module} we mean a module which is (irreducible, primitive)-$S$-atomic. Also we say two $S$-atomic
factorizations $m=s_1\cdots s_kn=t_1\cdots t_l n'$ are \emph{isomorphic}, if $k=l$, $n\sim^S n'$ and for a
permutation $\sigma$ of $\{1,\ldots,k\}$, we have $s_i\sim t_{\sigma(i)}$ for all $1\leq i\leq k$.

We say that $M$ is $S$-\pre\ when from $sm=m$ ($s\in S, m\in M$), we can deduce that $s\in \U(R)$ or $m=0$. By
\cite[Theorem 2.7(ii)]{smcs}, this is equivalent to saying that the three relations $\sim^S$, $\approx^S$ and
$\cong^S$ coincide or to asserting that $\cong^S$ is reflexive. In particular, all kinds of $S$-primitivity and
also by \cite[Theorem 2.7(iv)]{smcs}, all types of $S$-factorization are equivalent for a nonzero element of $M$,
if $M$ is $S$-\pre.

We call a module $M$, an \emph{$S$-unique factorization module}  or \emph{$S$-UFM} (resp. \emph{$S$-finite
factorization module} or \emph{$S$-FFM}), when $M$ is $S$-atomic and every nonzero element of $M$ has exactly one
(resp. finitely many) $S$-atomic factorization up to isomorphism. Also we say that $M$ is an \emph{$S$-bounded
factorization module} or \emph{$S$-BFM}, if for every $0\neq m\in M$ there is an $N_m\in \n$ such that the length
of every $S$-factorization of $m$ is at most $N_m$ and say that $M$ is an \emph{$S$-half factorial module} or
\emph{$S$-HFM} when $M$ is $S$-atomic and for each element $0\neq m\in M$ the length of all $S$-atomic
factorizations of $m$ are the same.

Note that in the cases that $S=R$ or $S=R=M$, these concepts coincide with the previously defined notations (see
\cite{r1,r2,we}). For example an integral domain $R$ is a UFD \ifof it is an $R$-UFM over itself. Moreover, a BFR
(bounded factorization Ring) means a ring which is a BFM over itself and a FFD (finite factorization domain) means
a domain $D$ which is a $D$-BFM. The notations UFR, FFR, BFD, HFD, \ldots, have similar meanings. Also, we
sometimes say \emph{$M$ has unique factorization (or has finite factorization or is \pre, \ldots) with respect to
$S$} instead of saying $M$ is an $S$-UFM (or $S$-FFM or $S$-\pre, \ldots).

Furthermore, if $E\se R$, we say that \emph{$R$ has unique factorization in $E$} when every nonzero nonunit
element in $E$ has unique factorization (with respect to $S=R$). Similar notations are used for other
factorization properties. In the following remark, we collect some observations which will be used in the paper
without any further mention.
\begin{rem}\label{easy rem}
\begin{enumerate}
\item \label{ER1} Every $S$-UFM is both an $S$-FFM and an $S$-HFM by definition and every $S$-BFM is $S$-\pre\
    (see remarks on page 8 of \cite{smcs}).
\item \label{ER2} If $R$ has unique factorization in $E$, then it is half factorial and has finite factorization
    in $E$ and if $R$ has bounded factorization in $E$, then it is \pre\ in $E$. Also if $E$ is a saturated
    multiplicatively closed subset or more generally, has the property that $xy\in E$ leads to $x\in E$ and
    $y\in E$, then being half factorial or having finite factorization in $E$ results to having bounded
    factorization in $E$ (see the second paragraph of \cite[p. 8]{smcs}).
\item \label{ER3} It is straightforward to see if $M=R$ and $s\in S$, then all kinds of $S$-primitivity for $s$
    are equivalent to being a unit and for elements in $S$, any type of $S$-associativity is equivalent to the
    corresponding type of $R$-associativity, since $S$ is assumed to be saturated.
\item \label{ER4} An element $m\in M$ is $S$-very strongly primitive, \ifof from $m=sm'$ for some $s\in S$ and
    $m'\in M$, we can deduce $s\in \U(R)$.
\end{enumerate}
\end{rem}

It should be mentioned that one can define other kinds of isomorphisms using different types of associativity and
also many forms of UFM, HFM, \ldots\ based on the choice of the type of irreducibility, primitivity and
isomorphism (see \cite{r1,r2} for the case $S=R$). But in order not to make the paper too long, we just
investigate the forms defined above, mentioning that similar techniques could be utilized to get similar results
on the other forms.

                   \section{$M$-splitting multiplicatively closed subsets}

A main concept used in \cite{d2} to relate factorization in $R$ and $R_S$ is the notion of a splitting multiplicatively closed subset of
$R$. A saturated multiplicatively closed subset $S$ of a domain $R$ is called a \emph{splitting multiplicative set}, when for each $x\in
R$, $x=as$ for some $a\in R$ and $s\in S$ such that $aR\cap tR=atR$ for all $t\in R$. An equivalent condition is
that principal ideals of $R_S$ contract to principal ideals in $R$ \cite[Lemma 1.2]{d2}. Here we will restate this
condition using factorization properties of the $R$-module $R$ with respect to $S$ and generalize it to every
$R$-module $M$. For this we need some more definitions.

\begin{defn}
By a compact $S$-atomic factorization of an element $m\in M$, we mean an equation of the form $m=sn$ for $s\in S$
and $S$-primitive element $n\in M$. We say that a subset $E\se M$ is compactly $S$-atomic if every nonzero element
of $E$ has a compact $S$-atomic factorization. If $E$ is compactly $S$-atomic and for every $0\neq m\in E$ and
compact $S$-atomic factorizations $m=sn=s'n'$ of $m$, we have $s\sim s'$ (resp. $s\sim s'$ and $n\sim^s n'$), then
$E$ is called semi-$S$-factorable (resp. $S$-factorable).
\end{defn}

Clearly every $S$-atomic module is compactly $S$-atomic but the following example shows that the converse is not
true. This example also shows that not all $S$-UFM's are factorable. Note that as usual when $S=R$, we drop the
$S$ prefixes.
\begin{ex}
Let $R$ be a ring with no irreducible elements (such as the domain $D$ in  \cite[Example 2.14]{smcs}) and $S=R$,
then the $R$-module $R$ is not atomic but is compactly atomic and even factorable, since $r=r1$ is the only
compact atomic factorization of an $r\in R$ up to associates. Also if $M=R/\fm$ for a maximal ideal $\fm$ of $R$,
then as in \cite[Example 2.14]{smcs}, $M$ is a UFM which is not even semi-factorable, since for any nonunit $r\in
R\sm \Ann(M)$ and any $m\in M\sm \Ann_M(r)$, $m'=rm=1m'$ are two compact atomic factorizations of $m'$ and $r\nsim
1$. Note that if we choose $R$ to be a valuation domain of Krull dimension 1 which is not a discrete valuation
domain, then $M$ is \pre.
\end{ex}
On the other hand, if $M$ is an $S$-UFM and $R$ is atomic in $S\sm \Ann(M)$, then $M$ is $S$-factorable. Because
if $m=sn$ is a compact $S$-atomic factorization of a nonzero element $m\in M$, then by replacing $s$ with its
atomic factorization, we get the unique $S$-atomic factorization of $m$ and hence $s$ and $n$ are unique up to
$S$-associates. Some properties of semi-$S$-factorable modules is stated in the next proposition.
\begin{prop}\label{semi elem}
Suppose that $E$ is a semi-$S$-factorable subset of $M$.
\begin{enumerate}
\item \label{2} If $m\in E$ is $S$-primitive, then $m\cong^S m$ and $m$ is $S$-very strongly primitive.
\item \label{3} If $S\cap \Z(M)=\tohi$, then $E$ is $S$-factorable.
\item \label{4} If $E=M$, then $R$ is \pre\ in $S\sm \Z(M)$, all kinds of irreducibility are equivalent for
    elements in $S\sm \Z(M)$ and all kinds of associativity are equivalent for pairs of elements in $S\sm
    \Z(M)$.
\end{enumerate}
\end{prop}
\begin{proof}
\ref{2}: It is easy to see that $m\cong^S m$. Now if $m=sm'$ for some $s\in S, m'\in M$, then since $m$ is
$S$-primitive, there is an $s'\in S$ such that $m'=s'm$, hence $m=ss'm$ and by $m\cong^S m$, it follows that
$s,s'\in \U(R)$ and $m$ is $S$-very strongly primitive.

\ref{3}: Let $0\neq m=sn=s'n'$ be two compact $S$-atomic factorizations of $m$. By semi-$S$-factorability, $s\sim
s'$ and hence $s=s''s'$ for some $s''\in R$. Since $S$ is saturated, $s''\in S$. Now $s's''n=s'n'$ and as
$s'\notin \Z(M)$, we get $s''n=n'$. Since $n'$ is $S$-very strongly primitive by \ref{2}, we deduce that $s''\in
\U(R)$ and hence $n'\sim^S n$, as required.

\ref{4}: Assume that $s=s's$ for some $s\in S\sm \Z(M), s'\in R$ and let $m$ be an $S$-primitive element of $M$.
Since $S$ is saturated, $s'\in S$ and because $sm=ss'm$  and $s\notin\Z(M)$, we get $m=s'm$ and by \ref{2}, it
follows that $s'\in \U(R)$. Therefore, $R$ is \pre\ in $S\sm \Z(M)$. Other parts of the claim follows from
\cite[Theorem 2.2(2)]{r1} or \cite[Theorem 2.7(iv)]{smcs}.
\end{proof}

Part \ref{3} of the above proposition shows that if $S\cap\Z(M)=\tohi$, then semi-$S$-factorability and
$S$-factorability are equivalent. Indeed, the author does not know an example of a semi-$S$-factorable module
which is not $S$-factorable even when $S\cap Z(M)$ is nonempty.

The next theorem and the remark following it, state conditions under which $S$-factorization properties of $M$ are
determined by factorization properties of elements in $S\sm \Ann(M)$.

\begin{thm}\label{S-X=X in S}
Suppose that $M$ is semi-$S$-factorable, $S\cap \Z(M)=\tohi$ and let $\cp$ be one of the following properties:
being \pre, having unique factorization, having finite factorization, being half factorial, having bounded
factorization, being atomic. Then $M$ has $\cp$  with respect to $S$ \ifof $R$ has $\cp$ in $S\sm \Ann(M)$.
\end{thm}
\begin{proof}
(\give): Let $s\in S\sm \Ann(M)$ and $m\in M\sm\Ann_M(s)$. By replacing $m$ with an $S$-primitive element
appearing in its compact $S$-atomic factorization, we can assume $m$ is $S$-primitive. First assume $\cp=$
atomicity. So $M$ is $S$-atomic and $sm$ has an $S$-atomic factorization $sm=s_1\cdots s_lm'$ with each $s_i$
irreducible and $m'$, $S$-primitive. Since $M$ is semi-$S$-factorable, we deduce that $s\sim s_1\cdots s_l$. By
\ref{semi elem}\ref{4}, $s\cong s_1\cdots s_l$, that is, $s=us_1\cdots s_l$ for some $u\in \U(R)$ and $s$ has an
atomic factorization.

Now for the other factorization properties, note that multiplication by $m$ turns any factorization of $s$ into an
$S$-factorization of $sm$ with the same length. Also this operation preserves isomorphism, so factorization
properties of $M$ pass to $S\sm \Ann(M)$. A similar argument takes care of $\cp=$ being \pre.

(\rgive): For $\cp$ = atomicity, the result is clear. Suppose that $R$ is \pre\ in $S\sm \Ann(M)$, $0\neq m=sm$
for some $m\in M$ and $s\in S$ and let $m=s'm'$ be a compact $S$-atomic factorization of $m$. Then $s'm'=ss'm'$
are two compact $S$-atomic factorizations for $m$ and hence $s'\sim ss'$. In particular, $s'=rss'$ for some $r\in
R$. Since $R$ is \pre\ in $S$, we deduce that $r,s\in \U(R)$, as required.

So assume $\cp\neq $ atomicity or being \pre. In any of the cases, $R$ is \pre\ in $S\sm \Ann(M)$ by \ref{easy
rem}\ref{ER2} and hence $M$ is $S$-\pre\ by the previous paragraph. Now if $x=s_1\cdots s_lm=s'_1\cdots s'_{l'}m'$
are two $S$-atomic factorizations of $0\neq x\in M$, then by semi-$S$-factorability of $M$, we get two $S$-atomic
factorizations of $s=s_1\cdots s_l\cong s'_1\cdots s'_{l'}$ with lengths $l, l'$. Thus if $\cp$= being half
factorial, then $l=l'$ and hence $M$ is an $S$-HFM. The case of bounded factorization is quite similar. For $\cp=$
having unique factorization or finite factorization, note that according to \ref{semi elem}\ref{3}, $m\sim^S m'$
in the above factorizations of $x$ and so if these two factorizations are non-isomorphic, then the two
factorizations of $s$ are also non-isomorphic. So the number of non-isomorphic factorizations of $x$ and $s$ are
the same.
\end{proof}

In several parts of the proof of the above result, we did not use the assumptions $S\cap \Z(M)=\tohi$ or
semi-$S$-factorability of $M$. So we get the following remark that states some weaker conditions under which, some
$S$-factorization properties of $M$ are determined by factorization properties of elements in $S\sm \Ann(M)$.
\begin{rem}\label{rem S-X=X in S}
Note that in the proof of \ref{S-X=X in S}(\give), for $\cp$= having bounded factorization or being \pre, we did
not use any of the two assumptions. Also for other properties, we could replace the condition ``$S\cap
\Z(M)=\tohi$'' with the weaker condition ``$R$ is \pre\ in $S\sm \Ann(M)$''.

In the proof of \ref{S-X=X in S}(\rgive), for $\cp$= being \pre, atomic or half factorial or having bounded
factorization, we did not need the condition ``$S\cap \Z(M)=\tohi$'' and for the other properties we could replace
the two conditions with ``$M$ is $S$-factorable.'' \qed
\end{rem}

Combining \ref{semi elem}\ref{4} with  the case $\cp=$ being \pre\ of \ref{S-X=X in S} we get:
\begin{cor}\label{split=> pre}
If $M$ is semi-$S$-factorable and $S\cap \Z(M)=\tohi$, then $M$ is $S$-\pre.
\end{cor}


Next we define the main concept of this research, namely $M$-splitting sets.
\begin{defn}\label{split def}
Let $E\se M$. We say that $S$ splits $E$ or $S$ is $E$-splitting, when the following two conditions hold.
\begin{enumerate}
\item \label{def1} $E$ is semi-$S$-factorable.

\item \label{def2} For every $S$-primitive element $r\in R$ and $S$-primitive element $m\in M$ such that $0\neq
    rm\in E$, the element $rm$ is $S$-primitive.
\end{enumerate}
\end{defn}

To see an example of this concept, see Section 5. The following result shows that this definition generalizes the
concept of splitting multiplicative sets as defined in \cite{d2}.
\begin{thm}\label{R-split=splitting}
Suppose that $R$ is a domain. Then $S$ is a splitting multiplicative set (in the sense of \cite{d2}) \ifof $S$
splits $R$.
\end{thm}
\begin{proof}
(\give): Note that since $R$ is a domain, it is \pre\ and all kinds of associativity are equivalent and also all
kinds of primitivity are equivalent. Suppose that $r\in R$ is $S$-primitive. By assumption we can write $r=sr'$
with $s\in S, r'\in R$ such that $Rr'\cap Rt=Rtr'$ for all $t\in S$. By $S$-primitivity, $r=ur'$ for some $u\in
\U(R)$. It follows that if $r\in R$ is $S$-primitive, then
\begin{align}
 Rr\cap Rt=Rtr\hbox{ for all }t\in S. \tag{$*$}
\end{align}
Conversely, if $r\in R$ satisfies $(*)$ and $r=sr'$ for some $s\in S, r'\in R$, then $r\in Rr\cap Rs=Rrs$, that
is, $r=r''rs$ for some $r''\in R$ and hence $1=r''s$, $s\in \U(R)$ and $r$ is $S$-primitive. Thus satisfying $(*)$
is equivalent to being $S$-primitive. Consequently, according to \cite[Corollary 1.4(a)]{d2}, $R$ is
$S$-factorable. Now assume that $r_1,r_2$ are nonzero $S$-primitive elements of $R$. By the above remarks $r_1$
and $r_2$ satisfy $(*)$ and hence by \cite[corollary 1.4(b)]{d2}, $r_1r_2$ also satisfies $(*)$ and hence is
$S$-primitive, as required.

(\rgive): It suffices to show that $S$-primitive elements of $R$, such as $r$, satisfy $(*)$. Let $0\neq x\in
Rr\cap Rt$, say $x=r_1r=r_2t$ for $r_1,r_2\in R$. If $r_i=s_ir'_i$ is the compact $S$-atomic factorization of
$r_i$ ($i=1, 2$), then $x=s_1(r'_1r)=(ts_2)r'_2$ are compact $S$-atomic factorizations of $x$, because by
assumption $r'_1r$ is $S$-primitive. So by semi-$S$-factorability, $s_1\sim ts_2$ and $s_1=t'ts_2$ for some $t'\in
R$. Thus $x=t'ts_2r'_1r\in Rtr$ and therefore, $Rr\cap Rt=Rtr$.
\end{proof}
At the end of this section we state a proposition which will be needed in the later sections.
\begin{prop}\label{M split -> R split}
Suppose that $S$ is $M$-splitting, $S'$ is a saturated multiplicatively closed subset of $R$ and $S'\sm \Ann(M)$
is compactly $S$-atomic. Then $S$ splits $S'\sm \Ann(M)$. If $S'=R$, then $S$ splits the $R$-module $R/\Ann(M)$.
\end{prop}
\begin{proof}
Suppose $r=s_1r_1=s_2r_2$ are two compact $S$-atomic factorizations of $r\in S'\sm\Ann(M)$. There exists an
$S$-primitive $m\in M\sm \Ann_M(r)$. Then $0\neq rm=s_1(r_1m)=s_2(r_2 m)$ are two compact $S$-atomic
factorizations of $rm$. So $s_1\sim s_2$. Now assume $r,r'\in S' \sm \Ann(M)$ are $S$-primitive and $rr'\in S'\sm
\Ann(M)$. If $rr'=sr''$ is a compact $S$-atomic factorization of $rr'$ and $m\in M\sm \Ann_M(rr')$ is
$S$-primitive, then $rr'm$ and $r''m$ are both $S$-primitive by condition \ref{def2} of definition of
$M$-splitting saturated multiplicatively closed subsets and hence from $rr'm=sr''m$ we deduce that $s\sim^S 1$ is
a unit. Therefore $rr'$ is $S$-primitive and $S$ splits $S'\sm \Ann(M)$.

To prove the claim about $N=R/\Ann(M)$, it suffices to show that an $r\in R\sm \Ann(M)$ is $S$-primitive \ifof its
image $\bar{r}$ is $S$-primitive in $N$. Assume that $r\in R\sm \Ann(M)$ is $S$-primitive and $\bar{r}=s\bar{r}'$.
So $r=sr'+a$ with $a\in \Ann(M)$. Let $m\in M$ be such that $rm\neq 0$. Since $M$ is compactly $S$-atomic, we can
assume that $m$ is $S$-primitive. Also suppose that $r'=s'r''$ is a compact $S$-atomic factorization of $r'$. Then
$rm=sr'm=ss'(r''m)$ and it follows from $S$ being $M$-splitting that $ss'\sim^S 1$, and $s\in \U(R)$. Thus
$\bar{r}$ is $S$-primitive. The reverse implication is straightforward.
\end{proof}

                              \section{Behavior of $S$-factorizations under localization}

Throughout this section we assume that $S\se S'$ are two saturated multiplicatively closed subsets of $R$ and set
$T$ to be the saturated multiplicatively closed subset, $S^{-1}S'= \{\f{s'}{s}| s'\in S', s\in S\}$ of $R_S$. We
investigate how factorization properties of $M$ with respect to $S'$ is related to factorization properties of
$M_S$ with respect to $T$, under the assumption that $S$ splits $M$. As we will see, in the case that $S'=R$, we
get some Nagata type theorems and also our results serve as partial answers to \cite[Question 3.11]{smcs}. To this
end, we first study how irreducibility behaves under localization.

\begin{prop}\label{irr}
Suppose that $S$  splits $E=S'\sm \Ann(M)$ and  $S\cap\Z(M)=S\cap \Z(R)=\tohi$. Let $\alpha\in \irr$ and $r=sa$ be
the compact $S$-atomic factorization of $r\in E\sm S$. Then $\bar{r}=r/1\in R_S$ is $\alpha$ \ifof $a$ is so in
$R$.
\end{prop}
\begin{proof}
Suppose that $a$ is very strongly irreducible and $\bar{r}=(r_1/s_1)(r_2/s_2)$. Note that as $r_1r_2=s_1s_2r\in
S'$, we must have $r_1,r_2\in S'$.  As $S\cap \Z(M)=\tohi$ and $r\notin \Ann(M)$, we conclude that $s_1s_2r\notin
\Ann(M)$ and hence $r_1,r_2\in E$. If $r_i=s'_ir'_i$ is the compact $S$-atomic factorization of $r_i$ ($i=1,2$),
then $s_1s_2sa=s'_1s'_2(r'_1r'_2)$ and thus $r'_1r'_2\in E$. So $r'_1r'_2$ is $S$-primitive by \ref{split
def}\ref{def2}. According to \ref{semi elem}\ref{3} (applied with $M=R$), we have $a\sim^Sr'_1r'_2$, in
particular, $a=tr'_1r'_2$ for some $t\in S$. Because $a$ is very strongly irreducible, one of $r'_i$'s, say $r'_1$
is a unit. Thus $r_1/s_1=s'_1r'_1/s_1\in \U(R_S)$ and as $\bar{r}$ is not a unit of $R_S$ (since $r\notin S$), it
is very strongly irreducible. Similar arguments show that if $a$ is (strongly) irreducible, then $\bar{r}$ is so.

Conversely, suppose that $\bar{r}$ is very strongly irreducible and $a=bc$. Then $b,c\in E$. So they have compact
$S$-atomic factorizations $b=s_1b'$ and $c=s_2c'$ with $b',c'\in E$. Since $S$ splits $E$ and $a=s_1s_2(b'c')$ is
a compact $S$-atomic factorization of the $S$-primitive element $a$ we deduce that $s_1,s_2\in \U(R)$ and hence
$b,c$ are both $S$-primitive. On the other hand, $\bar{r}=(\bar{s}\bar{b})\bar{c}$ and it follows from very
strongly irreducibility of $\bar{r}$ that for example $\bar{b}\in \U(R_S)$. This means that $b\in S$ and since $b$
is $S$-primitive, we conclude that $b\in \U(R)$. Therefore, $a$ is very strongly irreducible.

Now assume that $\bar{r}$ is strongly irreducible and $a=bc$. As in the above paragraph we see that $b,c$ are
$S$-primitive and $\bar{r}=(\bar{s}\bar{b})\bar{c}$. So by strongly irreducibility of $\bar{r}$ it follows that
for example $\bar{r}\approx \bar{b}$, that is, $\bar{r}=(s_1/s_2)\bar{b}$ for some $s_1,s_2\in S$. So $s_2sa=s_1b$
and by \ref{semi elem}\ref{3} $a\sim^S b$. But \ref{semi elem}\ref{2} implies that $a\cong^S b$ and hence
according to \cite[Theorem 2.7(i)]{smcs}, $a\approx b$, as required. The case $\alpha=$ irreducible is similar.
\end{proof}

Next we consider how $S'$-primitivity behaves under localization. Recall that throughout this section $S\se S'$
are saturated multiplicatively closed subsets and $T=S^{-1}S'$.
\begin{prop}\label{prim}
Suppose $S$ splits $M$, $S\cap \Z(M)=\tohi$ and $S'\sm \Ann(M)$ is compactly $S$-atomic. Let $0\neq m\in M$,
$\beta\in \prim$ and assume that $m=sn$ is the compact $S$-atomic factorization of $m$. Then $m/1$ is $T$-$\beta$
in the $R_S$-module $M_S$ \ifof $n$ is $S'$-$\beta$ in $M$.
\end{prop}
\begin{proof}
As $S\cap \Z(M)=\tohi$, we have $M\se M_S$, so we write $m$ instead of $m/1$. Assume that $m$ is $T$-$\beta$ and
$n=s'n'$ for $s'\in S', n'\in M$. If $n'=s_1n''$ and $s'=s_2s''$ are compact $S$-atomic factorizations of $n'$ and
$s'$, then we get $n=s_1s_2(s''n'')$ and as $S$ splits $M$, we deduce that $s_1s_2\sim 1$. This means that
$s_1,s_2\in \U(R)$ and both $n'$ and $s'$ are $S$-primitive. Also $m=ss'n'$.

Consider the case that $\beta=$ primitive. Then we get $m\sim ^T n'$, that is, $n'=(s'_0/s_0)m$ for $s'_0\in S'$,
$s_0\in S$. This leads to $s_0n'=s'_0sn=sv(s'_1n)$, where $s'_0=vs'_1$ is a compact $S$-atomic factorization of
$s'_0$ with $v\in S$. Consequently we get $n'\sim^S s'_1n$ and $n'=us'_1n$ for some $u\in S$. As $us'_1\in S'$, we
conclude that $n'\sim^{S'} n$ and hence $n$ is $S'$-primitive.

Now consider the case that $\beta=$ strongly primitive. Then we get $m\approx ^T n'$, whence $n'=(s'_0/s_0)m$ with
both $s_0,s'_0\in S$. Thus $s_0n=s'_0sn$ and as $S$ splits $M$, $n\sim^S n'$ which implies $n\cong^S n'$ by
\ref{semi elem}\ref{2} and hence by \cite[Theorem 2.7(i)]{smcs}, $n\approx^{S'} n'$, as required. We leave the
similar proof of the case $\beta=$ very strongly primitive to the reader.

Conversely, assume that $n$ is $S'$-primitive and $m=tm'$ for some $t\in T$ and $m'\in M_S$. One can readily check
that if $x=y/v$ for $y\in M, v\in S$ and $m\sim^T y$ (resp. $m\approx^T y$, $m\cong^T y$), then $m\sim^T x$ (resp.
$m\approx^T x$, $m\cong^T x$). Therefore, we can assume that $m'\in M$ and also $t=u'/u$ for $u'\in S'$ and $u\in
S$. If $u'=u_1u''$ and $m'=u_2m''$ are compact $S$-atomic factorizations of $u'$ and $m'$ with $u_1,u_2\in S$,
then it follows that $usn=um=u'm'=u_1u_2(u''m'')$ and hence $n\sim^S u''m''$. In particular, $n=s_0u''m''$ for
some $s_0\in S$. Note that $u''$ and hence $s_0u''$ are in $S'$, so as $n$ is $S'$-primitive, we deduce that
$n\sim^{S'} m''$, hence $m''=s'n$ for some $s'\in S'$. Thus $sm'=su_2m''=su_2s'n=u_2s'm$ and $m'=(u_2s'/s)m$.
Since $u_2s'/s\in T$, we see that $m\sim^T m'$ which shows that $m$ is $T$-primitive. The proof for (very)
strongly primitivity is similar.
\end{proof}


To see how $S'$-atomicity of $M$ and $T$-atomicity of $M_S$ are related, we need a lemma.
\begin{lem}\label{S'-prim => S-prim}
Suppose that $s'\in S'\sm(S\cup \Ann(M))$ and $0\neq m\in M$ such that $s'$ is $\alpha$ and $m$ is $S'$-$\beta$
where $\alpha\in \irr$, $\beta\in \prim$. If $M$ is semi-$S$-factorable, then $m$ is $S$-primitive. If $S$ splits
$M$ and $S'\sm \Ann(M)$ is compactly $S$-atomic, then $s'$ is $S$-primitive.
\end{lem}
\begin{proof}
If $m=sm'$ is a compact $S$-atomic factorization of $m$, then by $S'$-primitivity, $m\sim^{S'} m'$. So $m'=tm$ for
some $t\in S'$. Thus $m'=s(tm')$ and since $m'$ is $S$-very strongly primitive by \ref{semi elem}, we deduce that
$s$ is a unit and hence $m$ is $S$-primitive.

Now assume that $S$ splits $M$, $S'\sm \Ann(M)$ is compactly $S$-atomic and let $s'=s_1a$ be a compact $S$-atomic
factorization of $s'$. Since $s'$ is irreducible, either $s'\sim s_1$ or $s'\sim a$. In the former case, it
follows that $s_1=rs'$ and as $S$ is saturated, we get the contradiction $s'\in S$. So $s'\sim a$ and $a=rs'$ for
some $r\in R$. Note that since $s'=s_1s'r$ and $S'$ is saturated, we must have $r\in S'$. If $r=s_2r'$ is a
compact $S$-atomic factorization of $r$, then $a=s_1s_2(r'a)$ and thus $s_1s_2\sim 1$, for both sides are compact
$S$-atomic factorizations in $E=S'\sm \Ann(M)$ and $S$ splits $E$ by \ref{M split -> R split}. This means that
$s_1\in \U(R)$ and the result follows.
\end{proof}
\begin{thm}\label{atomic}
Suppose that $S$ is $M$-splitting, $S'\sm\Ann(M)$ is compactly $S$-atomic, $\alpha\in \irr$ and $\beta\in \prim$.
Then the following hold.
\begin{enumerate}
\item \label{atomic1} If $M$ is ($\alpha$, $\beta$)-$S'$-atomic, then $M$ is ($\alpha$, very strongly
    primitive)-$S$-atomic.

\item \label{atomic3} Assume that $S \cap \Z(M)=S\cap\Z(R)=\tohi$. Then $M$ is ($\alpha$, $\beta$)-$S'$-atomic
    \ifof $M$ is ($\alpha$, primitive)-$S$-atomic and $M_S$ is ($\alpha$, $\beta$)-$T$-atomic.
\end{enumerate}
\end{thm}
\begin{proof}
\ref{atomic1}: If $m=s_1\cdots s_ks'_1\cdots s'_{k'}m'$ is an ($\alpha$, $\beta$)-$S'$-atomic factorization of
$0\neq m\in M$ where $s_i\in S$ and $s'_i\in S'\sm S$, then by \ref{S'-prim => S-prim}, $s'_i$ and $m'$ are
$S$-primitive and hence their product $0\neq m''=s'_1\cdots s'_{k'}m'$ is also $S$-primitive and by \ref{semi
elem}\ref{2}, indeed $S$-very strongly primitive. Therefore $m=s_1\cdots s_k m''$ is an ($\alpha$, very strongly
primitive)-$S$-atomic factorization of $m$.

\ref{atomic3}: (\give): It follows from \ref{atomic1} that   $M$ is ($\alpha$, primitive)-$S$-atomic. Suppose
$m=s_1\cdots s_ks'_1\cdots s'_{k'}m'$ is an ($\alpha$, $\beta$)-$S'$-atomic factorization of $0\neq m\in M$ where
$s_i\in S$ and $s'_i\in S'\sm S$, then by \ref{S'-prim => S-prim}, \ref{prim} and \ref{irr}, $s'_i/1$ is $\alpha$
and $m'/1$ is $T$-$\beta$. Hence we get the ($\alpha$, $\beta$)-$T$-atomic factorization $m/1= (s'_1/1) \cdots
(s'_{k'}/1) (um'/1)$, where $u=(s_1\cdots s_k)/1\in \U(R_S)$. So  $M_S$ is ($\alpha$, $\beta$)-$T$-atomic.

(\rgive): Let $0\neq m\in M$. Then by assumption $m=s_1\cdots s_k n$ where $s_i\in S$ are $\alpha$ and $n$ is
$S$-primitive. Assume $n/1=(s'_1/u_1)\cdots (s'_{k'}/u_{k'}) (m'/u_{k'+1})$ is an ($\alpha$, $\beta$)-T-atomic
factorization where $s'_i\in S'\sm S$, $u_i\in S$ and let $s'_i=v_is''_i$ and $m'=vm''$ be compact $S$-atomic
factorizations of $s'_i$ and $m'$, respectively. Then as $s'_i/u_i$  and hence $s'_i/1$ are $\alpha$ it follows
from \ref{irr} that $s''_i$ is $\alpha$. Similarly by \ref{prim}, $m''$ is $S'$-$\beta$. Since all $s''_i$ and
$m''$ are $S$-primitive and $S$ splits $M$, $n'=s''_1\cdots s''_{k'}m''$ is $S$-primitive. Now from the above
factorization of $n/1$ it follows that $un=u'n'$ for some $u,u'\in S$ and since $n$ and $n'$ are both
$S$-primitive and by \ref{semi elem}\ref{3}, $n\sim^S n'$ and hence $n\cong ^S n'$ according to \ref{semi
elem}\ref{2}. This implies that $n=an'$ with $a\in \U(R)$. We conclude that $m= s_1\cdots s_k s''_1\cdots s''_{k'}
(am'')$ is an ($\alpha$, $\beta$)-$S'$-atomic  factorization of $m$.
\end{proof}


To establish a version of the above theorem for other factorization properties we need:
\begin{lem}\label{S'-HF FF=>S-BF}
Suppose that $M$ is semi-$S$-factorable and $S\cap \Z(M)=\tohi$. If $M$ is an $S'$-UFM, $S'$-HFM or $S'$-FFM, then
$M$ is an $S$-BFM. Also if $M$ is an $S$-BFM, then $E=R\sm \Ann(M)$ is compactly $S$-atomic.
\end{lem}
\begin{proof}
Let $0\neq m\in M$. In either of the cases, the possible lengths of an $S'$-atomic factorization of $m$ are finite
and hence there is an upper bound $N_m$ on these lengths. Now let $m=s_1\cdots s_lm'$ where $s_i\in S\sm\U(R)$. By
replacing $m'$ with one of its compact $S$-atomic factorizations, the length of this factorization does not
decrease. Therefore, we can assume that $m'$ is $S$-primitive.

Let $s_1m'=s'_{1,1}\cdots s'_{1,k_1}m'_1$ be an $S'$-atomic factorization of $s_1m'$. If $k_1=0$, then as $m'_1$
is $S$-very strongly primitive by \ref{S'-prim => S-prim} and \ref{semi elem}\ref{2}, we must have $s_1\in \U(R)$
against our assumption. So $k_1\geq 1$. Similarly, we can find $S'$-atomic factorizations
$s_im'_{i-1}=s'_{i,1}\cdots s'_{i,k_i}m'_i$ for each $2\leq i\leq l$ with $k_i\geq 1$ for all $i$. Consequently,
we get an $S'$-atomic factorization $m=s'_{1,1}\cdots s'_{1,k_1}\cdots s'_{l,1}\cdots s'_{l,k_l}m'_l$. Hence
$l\leq \sum_{i=1}^l k_i\leq N_m$ and $M$ is an $S$-BFM.

Now assume that $r\in E$ and $M$ is an $S$-BFM. Then there is an $S$-primitive element $m\in M\sm \Ann_M(r)$. If
$r=s_1\cdots s_l r'$ with $s_i\in S\sm \U(R)$, then $rm=s_1\cdots s_l (r'm)$ and hence $l\leq N_{rm}$. If $l$ is
the largest possible length of $S$-factorizations of $r$, then $r'$ is $S$-very strongly primitive and the result
follows.
\end{proof}
Note that although under the conditions of the first part of the above lemma, $M$ is an $S$-BFM, but it need not
be an $S'$-BFM as demonstrated in \cite[Example 2.14]{smcs} with $S'=D$ and $S=\U(D)$.

\begin{thm}\label{S'-P=>S-P}
Suppose that $M$ is semi-$S$-factorable, $S\cap\Z(M)=\tohi$ and let $\cp$ be one of the following properties:
being \pre, having unique factorization, having finite factorization, being half factorial, having bounded
factorization. If $M$ has $\cp$ with respect to $S'$ then it has $\cp$ with respect to $S$.
\end{thm}
\begin{proof}
For the case that $\cp=$ being \pre\ or having bounded factorization, the result is \cite[Theorem 3.8]{smcs} (and
does not need semi-$S$-factorability or $M$-regularity). For other cases, note that by the previous lemma, $M$ is
an $S$-BFM and hence $S$-atomic and therefore $R$ is atomic in $E=S'\sm\Ann(M)$ by \ref{S-X=X in S}. Let $s\in E$
and $m\in M\sm \Ann_M(s)$. By replacing $m$ with an $S'$-primitive element appearing in an $S'$-atomic
factorization of $m$, we can assume that $m$ is $S'$-primitive. Now if $s=s_1\cdots s_k$ is an atomic
factorization, then $sm=s_1\cdots s_km$ is an $S'$-atomic factorization of $sm$ with the same length and two
factorizations of $sm$ arising in this way are isomorphic \ifof the two factorizations of $s$ are isomorphic.
Consequently, $R$ has $\cp$ in $E$ and hence by \ref{S-X=X in S}, $M$ has $\cp$ with respect to $S$.
\end{proof}

Another condition under which, an $S'$-UFM is an $S$-UFM is presented in \cite[Theorem 3.8 \& Notation 3.5]{smcs}.
\begin{lem}\label{iso}
Assume that $S$ splits $M$, $S'\sm \Ann(M)$ is compactly $S$-atomic and $S\cap \Z(M)=S\cap\Z(R)=\tohi$. Let
$y_1,y_2\in S'\sm(S\cup \Ann(M))$ be irreducible and $0\neq m_1,m_2\in M$ be $S'$-primitive. Then $y_1/1 \sim
y_2/1$ in $R_S$ \ifof $y_1\sim y_2$ in $R$ and $m_1/1\sim ^T m_2/1$ \ifof $m_1\sim ^{S'} m_2$.
\end{lem}
\begin{proof}
(\rgive): Trivial. (\give): Suppose $y_1/1=(r/s) (y_2/1)$. Then $sy_1= ry_2=s_0(r'y_2)$, where $r=s_0r'$ is a
compact $S$-atomic factorization of $r$. Hence by \ref{S'-prim => S-prim}, \ref{M split -> R split} and \ref{semi
elem}\ref{3}, $y_1\sim^S r'y_2$ and $y_1\in Ry_2$. Similarly $y_2\in Ry_1$ and $y_1\sim y_2$. The proof of the
other statement is similar.
\end{proof}

\begin{thm}\label{S-P & S-1 P}
Suppose that $S$ splits $M$, $S\cap \Z(M)=S\cap \Z(R)=\tohi$ and let $\cp$ be one of the following properties:
having unique factorization, having finite factorization, being half factorial, having bounded factorization. Then
$M$ has $\cp$ with respect to $S'$ \ifof $M$ has $\cp$ with respect to $S$ and $M_S$ has $\cp$ with respect to
$T$. A similar statement holds for $\cp=$ being \pre, if we further assume that $S'\sm \Ann(M)$ is compactly
$S$-atomic.
\end{thm}
\begin{proof}
(\give): According to \ref{S'-P=>S-P}, we just need to show that $M_S$ has $\cp$ with respect to $T$. Note that by
\ref{S'-HF FF=>S-BF}, in all cases $S'\sm \Ann(M)$ is compactly $S$-atomic. First assume that $\cp=$ being \pre\
and $m/s_0= (s'/s_1) (m/s_0)$ for some $s'\in S', s_1,s_0\in S,0\neq m\in M$. If $m=v_0m'$ and $s'=v_1s''$ are
compact $S$-atomic factorizations of $m$ and $s'$ with $v_i\in S$, then $s_1m'=v_1(s''m')$ and as $S$ splits $M$
and by \ref{semi elem}, it follows that $m'\cong^{S} s''m'$. Then $m'=(us'')m'$ for some $u\in \U(R)$ and since
$M$ is $S'$-\pre\ and $s''\in S'$, we deduce $s''\in \U(R)$ and $s'\in S$, as required.

For $\cp\neq$ being \pre, since $M$ is $S'$-atomic, it follows form \ref{atomic} that $M_S$ is $T$-atomic. Let
$0\neq m\in M$ be $S$-primitive. If $(m/1)=(s'_1/s_1)\cdots (s'_k/s_k)(m'/1)$ is any $T$-factorization of $m/1$
and $s'_i= v_is''_i$ and $m'=vm''$ are compact $S$-atomic factorization of $s'_i$ and $m'$, then $s_1\cdots s_k
m=v_1\cdots v_k (s''_1\cdots s''_km'')$ which by \ref{semi elem}\ref{3} implies that $m\cong^{S} s''_1\cdots
s''_km''$. Also according to \ref{irr} and \ref{prim}, $s''_i$ is irreducible \ifof $s'_i/s_i$ is so and $m''$ is
$S'$-primitive \ifof $m'/1$ is so. Therefore, if for example the number of $S'$-atomic factorizations of $m$ is
finite, then so is the number of $T$-atomic factorizations of $m/1$. Similarly other $S'$-factorization properties
of $m$ pass to $T$-factorization properties of $m/1$. Noting that each $0\neq x\in M_S$ is a unit multiple of some
$m/1$ where $0\neq m\in M$ is $S$-primitive, the result is concluded.

(\rgive): The cases $\cp=$ being \pre\ or having bounded factorization is \cite[Theorem 3.9]{smcs} (with much less
assumptions). So assume that $\cp=$ having unique factorization or finite factorization or being half factorial.
Note that by \ref{S'-HF FF=>S-BF} applied with $S'=S$, we see that $M$ is $S$-BFM and $R\sm \Ann(M)$ is compact
$S$-atomic. Thus it follows \ref{atomic}\ref{atomic3}, that $M$ is $S'$-atomic. We prove the result for $\cp=$
having finite factorization and the other cases follow similarly.

Let $m=s_1\cdots s_k s'_1\cdots s'_{k'} m'$ be an $S'$-atomic factorization of $0\neq m\in M$ with $s_i\in S$ and
$s'_i\in S'\sm S$. Then by \ref{S'-prim => S-prim}, each $s'_i$, $m'$ and hence their product $m''=s'_1\cdots
s'_{k'}m'$ are $S$-primitive. Consequently, $m=s_1\cdots s_k m''$ is $S$-isomorphic (and thus $S'$-isomorphic) to
one of the finite $S$-atomic factorizations of $m$. So if we show $S$-primitive elements of $M$ have finitely many
$S'$-atomic factorizations, we are done. Thus we assume $m$ is $S$-primitive and $k=0$. Then $m/1=(s'_1/1) \cdots
(s'_{k'}/1) (m'/1)$ is a $T$-atomic factorization of $m/1$ by \ref{irr} and \ref{prim}. Therefore, each
$S'$-atomic factorization of $m$ leads to a $T$-atomic factorization of $m/1$ and according to \ref{iso} if two
such $T$-atomic factorizations of $m/1$ are $T$-isomorphic, then the original $S'$-atomic factorizations of $m$
are $S'$-isomorphic. Consequently, as $m/1$ has finitely many $T$-atomic factorizations up to $T$-isomorphisms,
$m$ also has finitely many $S'$-atomic factorizations up to $S'$-isomorphism, as claimed.
\end{proof}

It should be mentioned that (\rgive) of the above theorem, is a partial answer to \cite[Question 3.11]{smcs}.
Summing up Theorems \ref{atomic}\ref{atomic3}, \ref{S-P & S-1 P} and \ref{S-X=X in S}, we get:
\begin{cor}\label{main result}
Let $S\se S'$ be two saturated multiplicatively closed subsets of $R$ and set $T=S^{-1}S'$. Suppose that $S$
splits $M$ and $S\cap \Z(M)= S\cap \Z(R)=\tohi$. For $\cp\in$ \{having unique factorization, having finite
factorization, being half factorial, having bounded factorization\}, the following are equivalent.
\begin{enumerate}
\item $M$ has $\cp$ with respect to $S'$.

\item $M$ has $\cp$ with respect to $S$ and $M_S$ has $\cp$ with respect to $T$.

\item $R$ has $\cp$ in $S\sm \Ann(M)$ and $M_S$ has $\cp$ with respect to $T$.
\end{enumerate}
If we further assume that $S'\sm \Ann(M)$ is compactly $S$-atomic, then the above conditions are also equivalent
for $\cp=$ being \pre\ or atomic.
\end{cor}
If we set $S'=R$ in the above corollary, we get some Nagata type theorems. If we further assume that $M=R$ and $R$
is an integral domain, then this corollary implies \cite[Theorems 3.1 \& 3.3]{d2} (except for ACCP and idf-domain,
which are not investigated in this research). Because if $S$ is generated by primes, then every element of $S$ has
a unique factorization as a product of primes and hence irreducibles, that is, $R$ has unique factorization (and
whence has finite and bounded factorization and is \pre, atomic and half factorial) in $S$. Indeed, even when
$M=S'=R$ and $R$ is a domain, this corollary is slightly stronger than \cite[Theorem 3.1]{d2}, since in our
results $S$ need not be generated by primes. An example in which $S$ is not generated by primes is presented in
the next section.

                              \section{An application}

We present an example which shows how our main result \ref{main result}, could be applied in the case that
$M=S'=R$ and $R$ is an integral domain, the classical and most important situation in the factorization theory.
Note that in this case, $R$ is \pre\ and hence all types of associativity (resp. irreducibility, primitivity) are
equivalent to each other. Also if $S$ splits $R$, then $R$ is compactly $S$-atomic and hence \ref{main result}
could be applied for all $\cp\in$ \{having unique factorization, having finite factorization, being half
factorial, having bounded factorization, being atomic\}. First let's state the setting of our example as a
notation.
\begin{nota}\label{Ex nota}
In this section, we assume that $M=S'=R=A+XB[X]$ where $A\se B$ are integral domains and $S=(\U(B)\cap A)\cup
\{uX^n|u\in \U(B), n\in \n\}$. Also we set $S_0=\U(B)\cap A$ and denote the quotient field of $A$ by $K$.
\end{nota}

It is easy to see that $S$ is a saturated multiplicatively closed subset of $R$, indeed, it is the saturated
multiplicatively closed subset generated by $X$. Also if $B\neq A$ and $b\in B\sm A$, then $(bX)^2\in RX$, while
$bX\notin RX$. So $X$ is not prime and $S$ is not generated by primes if $A\neq B$.
\begin{thm}\label{Ex splits}
The set $S$ splits $R$ \ifof all of the following conditions hold:
\begin{enumerate}
\item \label{Ex1} $S_0$ splits $A$;
\item \label{Ex2} for every $b\in B$ there are $u\in \U(B)$ and $a\in A$ such that $b=ua$;
\item \label{Ex3} $\U(B)\cap K\se A_{S_0}$.
\end{enumerate}
In particular, if either $B$ is a filed or $S_0$ is any $A$-splitting saturated multiplicatively closed subset of
$A$ and $B=A_{S_0}$, then $S$ splits $R$.
\end{thm}
\begin{proof}
(\give): Note that if $a\in A$, and $a=sf$ for some $s\in S, f\in R$, then $s\in S_0$ and $f\in A$. Thus $a$ is
$S$-primitive \ifof it is $S_0$-primitive and so \ref{Ex1} follows as $S$ splits $A$. To see \ref{Ex2}, let $b\in
B$. If $bX$ is $S$-primitive, then $b^2X^2$ should be $S$-primitive by \ref{def2} of Definition \ref{split def}.
But $b^2X^2=(b^2X)X$ is an $S$-factorization of $b^2X^2$, a contradiction. So $bX$ is not $S$-primitive. Thus if
its compact $S$-atomic factorization  is $bX=sf$, with $s\in S$ and $S$-primitive $f\in R$, then $\deg f=0$.
Whence $f\in A$ and $s=uX$ for some $u\in B$ and \ref{Ex2} follows with $a=f$.

Now let $u\in \U(B)$ and $u=a/a'$ for some $0\neq a,a'\in A$. Assume that $a=s_0 a_0$ and $a'=s'_0a'_0$ are
compact $S$-atomic (hence $S_0$-atomic) factorizations of $a$ and $a'$. Then $((us'_0/s_0)X)a'_0=(X)a_0$ are two
compact $S$-atomic factorizations of $a_0X$, therefore by uniqueness of such factorizations we have $X\cong
(us'_0/s_0) X$, that is, $v= us'_0/s_0\in \U(R)=\U(A)$. Consequently, $u=vs_0/s'_0\in A_{S_0}$.

(\rgive): Suppose that $f=\sum_{i=k}^n b_i X^i\in R$ with $k\leq n$ and $b_k\neq 0$. If $k=0$, then $b_0\in A$ has
a compact $S_0$-atomic factorization $b_0=sa$ where $s\in \U(B)\cap A$. Therefore $f= s (a+\sum_{i=2}^n
(b_i/s)X^i)$ is a compact $S$-atomic factorization of $f$ (note that $a+\sum_{i=2}^n (b_i/s)X^i$ is $S$-primitive,
since $a$ is so). If $k>0$, then by \ref{Ex2}, there are $u\in \U(B)$ and $a_0\in A$ such that $b=ua_0$. If
$a_0=s_0 a$ is the compact $S_0$-atomic factorization of $a_0$, then set $s= us_0\in \U(B)$. Thus we get the
compact $S$-factorization $f= (sX^k) (a+\sum_{i=k+1}^n (b_i/s) X^{i-k})$. Thus $R$ is compactly $S$-atomic and $f$
is $S$-primitive \ifof $b_0=f(0)$ is nonzero and $S_0$-primitive in $A$. In particular, if $f, g\in R$ are
$S$-primitive, then $0\neq (fg)(0)$ is $S_0$-primitive in $A$ by \ref{Ex1} and $fg$ is $S$-primitive.

It remains to show the uniqueness of compact $S$-factorizations of $f$. Assume that $f=u_1X^{k_1} f_1$ and $u_2
X^{k_2} f_2$ are two $S$-atomic factorizations of $f$ where $f_1$ and $f_2$ are $S$-primitive and $u_1,u_2\in
\U(B)$. Therefore, by the previous paragraph, $a_1=f_1(0)\neq 0$, $a_2=f_2(0)\neq 0$ are $S_0$-primitive and hence
$k_1=k_2=k$. Also $u_1a_1=u_2a_2$ and hence $u= u_1/u_2= a_2/a_1\in \U(B)\cap K\se A_{S_0}$, by \ref{Ex3}.
Consequently, $u=a/s_0$ for some $a\in A$ and $s_0\in S_0$. If $a=s'_0a'$ is the compact $S_0$-atomic
factorization of $a$, then $s_0a_2= a a_1= s'_0(a'a_1)$. Since both $a_2$ and $a'a_1$ are $S_0$-primitive, we must
have $a_2\sim^{S_0} a'a_1$ and in particular, $a_2\in Aa_1$. Similarly $a_1\in Aa_2$, hence $Aa_1=Aa_2$ and as $A$
is a domain, $a_1=va_2$ for some $v\in \U(A)$. It follows that $u_2=vu_1$ and hence $u_1X^k\sim^S u_2X^k$, as
required.
\end{proof}

Thus in particular, we can apply \ref{main result} with $M=S'=R$, in the case $A=B$. Since in this case $X$ is
prime, $R=B[X]$ has unique factorization in $S$ and hence we get
\begin{cor}\label{B[X,X^-1]}
If $B$ is an integral domain, then $B[X]$ is atomic (resp. a BFD, a FFD, a HFD, a UFD) \ifof $B[X, X^{-1}]$ is so.
\end{cor}
\begin{thm}\label{Ex main}
Using Notation \oldref{Ex nota} and assuming that the conditions \ref{Ex1}--\ref{Ex3} of \ref{Ex splits} hold,
then we have
\begin{enumerate}
\item \label{Ex main 1} $R$ is atomic (resp. a BFD, a HFD) \ifof $U(B)\cap A= \U(A)$ and $B[X]$ is atomic (resp.
    a BFD, a HFD).

\item \label{Ex main 2} $R$ is a FFD \ifof $|\f{\U(B)}{\U(A)}|<\infty$ and $B[X]$ is a FFD.

\item \label{Ex main 3} $R$ is a UFD \ifof $|\f{\U(B)}{\U(A)}|=1$ and $B[X]$ is a UFD.
\end{enumerate}
\end{thm}
\begin{proof}
In all cases, by applying \ref{main result}, we deduce that $R$ has the desired property $\cp$ \ifof $R$ has $\cp$
in $S$ and $R_S=B[X,X^{-1}]$ has $\cp$. Thus according to \ref{B[X,X^-1]}, we just need to show that in each case,
$R$ has $\cp$ in $S$  \ifof the stated condition on $\U(B)$ and $\U(A)$ is satisfied.

\ref{Ex main 1}: Assume $R$ is atomic in $S$. Then $X$  has decomposition $X=a_1\cdots a_n (uX)$ with irreducible
$a_i\in A$ and $u\in \U(B)$ such that $uX$ is irreducible in $R$. If $a\in \U(B)\cap A$, then $uX=((u/a)X)a$ and
as $uX$ is irreducible, it follows that $a\in\U(R)=\U(A)$. Conversely, if $\U(B)\cap A=\U(A)$, then $uX$ is
irreducible for each $u\in \U(B)$ and hence $f= uX^n= (uX)X^{n-1}$ is an atomic factorization of $f\in S\sm\U(R)$.
Also in any atomic factorization of $f$ exactly $n$ terms of the form $vX$ appear where $v\in \U(B)$. Since
irreducible elements of $S$ are exactly those of the form $vX$ for $v\in \U(B)$, it follows that any atomic
factorization of $f$ has length $n$ and $R$ is half factorial and has bounded factorization in $S$.

\ref{Ex main 2}: If $R$ has finite factorization in $S$, then it is atomic in $S$ and hence by \ref{Ex main 1},
$\U(B)\cap A=\U(A)$ and the irreducible elements in $S$ are of the form $vX$ for $v\in \U(B)$. Note that $vX\cong
v'X$ \ifof $v\in \U(A)v'$ \ifof the image of $v,v'$ in the quotient group $\U(B)/\U(A)$ are equal. So if $v_1,
v_2, \cdots$ are infinite elements of $\U(B)$, such that $v_i\U(A)\neq v_j\U(A)$ for each $i\neq j$, then we can
find an infinite set of non-isomorphic atomic factorizations of $X^2=(v_1X)(v_1^{-1}X)= (v_2X)(v_2^{-1}X)=\cdots$.

Conversely, suppose that $n= |\U(B)/\U(A)|<\infty$. If $a\in\U(B)\cap A$, then $a^n\in \U(A)$  and hence $a\in
\U(A)$. Therefore, by \ref{Ex main 1} $R$ is atomic and half factorial in $S$ and irreducibles of $R$ in $S$ are
of the form $vX$ with $v\in\U(B)$. Since every atomic factorization of $uX^k\in S$ has the same length $k$, to
show that it has finitely many factorizations it suffices to show that it has only finitely many  non-associate
irreducible divisors. But there are $n$ non-associate irreducible elements in $S$, because $vX\cong v'X$ \ifof
$v\U(A)=v'\U(A)$, and we are done. The proof of \ref{Ex main 3} is similar.
\end{proof}
This theorem generalizes the following previously known result (for example, it is an immediate consequence of the
propositions considering the $D+M$ construction in \cite{d1}) which follows from \ref{Ex main} in the case that
$B$ is a field.
\begin{cor}
Assume that $B$ is a field. Then $R$ is atomic \ifof $R$ is a BFD \ifof $R$ is a HFD \ifof $A$ is a field. Also
$R$ is a FFD (resp. UFD) \ifof $A$ is a field and $|B^*/A^*|<\infty$ (resp. $B=A$).
\end{cor}

\subsection*{Acknowledgement} This research was financially supported by a grant of National Elites Foundation of
Iran. The author would like to thank the referee for his/her nice comments.
                                             
\end{document}